# LARGE DEVIATION ASYMPTOTICS AND CONTROL VARIATES FOR SIMULATING LARGE FUNCTIONS[1]

By Sean P. Meyn

*University of Illinois at Urbana–Champaign*

Consider the normalized partial sums of a real-valued function $F$ of a Markov chain,

$$\phi_n := n^{-1} \sum_{k=0}^{n-1} F(\Phi(k)), \qquad n \geq 1.$$

The chain $\{\Phi(k) : k \geq 0\}$ takes values in a general state space $\mathsf{X}$, with transition kernel $P$, and it is assumed that the Lyapunov drift condition holds: $PV \leq V - W + b\mathbb{I}_C$ where $V : \mathsf{X} \to (0, \infty)$, $W : \mathsf{X} \to [1, \infty)$, the set $C$ is small and $W$ dominates $F$. Under these assumptions, the following conclusions are obtained:

1. It is known that this drift condition is equivalent to the existence of a unique invariant distribution $\pi$ satisfying $\pi(W) < \infty$, and the law of large numbers holds for any function $F$ dominated by $W$:

$$\phi_n \to \phi := \pi(F), \qquad \text{a.s., } n \to \infty.$$

2. The lower error probability defined by $\mathsf{P}\{\phi_n \leq c\}$, for $c < \phi$, $n \geq 1$, satisfies a large deviation limit theorem when the function $F$ satisfies a monotonicity condition. Under additional minor conditions an exact large deviations expansion is obtained.

3. If $W$ is near-monotone, then control-variates are constructed based on the Lyapunov function $V$, providing a *pair of estimators* that together satisfy nontrivial large asymptotics for the lower *and* upper error probabilities.

In an application to simulation of queues it is shown that exact large deviation asymptotics are possible even when the estimator does not satisfy a central limit theorem.

Received October 2004; revised July 2005.

[1]Supported by the National Science Foundation under Awards ECS 02 17836 and DMI-00-85165.

**1. Introduction.** This paper explores extensions of the control-variate method to obtain confidence bounds in simulation of a function of a Markov chain $\boldsymbol{\Phi} = \{\Phi(0), \Phi(1), \ldots\}$. It is assumed that $\boldsymbol{\Phi}$ evolves on a general state space $\mathsf{X}$, equipped with a countably generated sigma-field $\mathcal{B}$. The statistics of $\boldsymbol{\Phi}$ are determined by its initial distribution, and the transition kernel $P$ defined by

$$P(x, A) := \mathsf{P}\{\Phi(1) \in A | \Phi(0) = x\}, \qquad x \in \mathsf{X}, A \in \mathcal{B}.$$

Let $\{L_n : n \geq 1\}$ denote the sequence of empirical measures induced by $\boldsymbol{\Phi}$ on $(\mathsf{X}, \mathcal{B})$,

$$(1) \qquad L_n := \frac{1}{n} \sum_{k=0}^{n-1} \delta_{\Phi(k)}, \qquad n \geq 1.$$

It is assumed that $\boldsymbol{\Phi}$ is positive Harris recurrent, with unique invariant probability distribution denoted $\pi$. Equivalently, for each bounded measurable function $F: \mathsf{X} \to \mathbb{R}$, and each initial condition, the law of large numbers holds:

$$L_n(F) \to \phi := \pi(F) = \int F(x) \pi(dx), \qquad n \to \infty, \text{ a.s.}$$

See [19, 22] or [38], Theorem 17.1.7. For each measurable function $F: \mathsf{X} \to \mathbb{R}$ satisfying $\pi(|F|) < \infty$, the sequence $\{L_n(F) : n \geq 1\}$ is interpreted as Monte Carlo estimates of the steady-state mean of $F$.

While consistent for each initial condition even when $F$ is not bounded, finer assumptions are required to obtain *confidence bounds*, that is, bounds on $\mathsf{P}\{|L_n(F) - \phi| \geq a\}$ for a given $a > 0$. Such bounds are typically based on one of the following limit theorems:

THE CENTRAL LIMIT THEOREM (CLT). *For some $\sigma \geq 0$ and each initial condition,*

$$(2) \qquad \sqrt{n}[L_n(F) - \phi] \xrightarrow{w} \sigma X, \qquad n \to \infty,$$

*where $X$ is a standard normal random variable, and the convergence is in distribution.*

THE LARGE DEVIATION PRINCIPLE, OR LDP. *For a convex function $I : \mathbb{R} \to \mathbb{R}_+$, and any nonempty open interval $(c_0, c_1) \subset \mathbb{R}$,*

$$
\begin{aligned}
&\lim_{n \to \infty} n^{-1} \log \mathsf{P}\{L_n(F) - \phi \in (c_0, c_1)\} \\
(3) \qquad &= \lim_{n \to \infty} n^{-1} \log \mathsf{P}\{L_n(F) - \phi \in [c_0, c_1]\} \\
&= -\left(\min_{c \in [c_0, c_1]} I(c)\right).
\end{aligned}
$$

There has been tremendous research activity concerning large deviation properties of Markov chains following the pioneering work of Donsker and Varadhan [15, 46, 47]. The literature contains a broad range of possible conclusions under a correspondingly broad range of assumptions (see the monographs [11, 12]).

The strongest conclusions are based on variants of the assumptions imposed by Donsker and Varadhan in [13, 14, 15], that are essentially equivalent to *compactness* of the $n$-step transition operator, for some $n > 0$ (see [45], Theorem 2.1, or [10], Lemma 3.4). Under these assumptions the LDP holds for the empirical distributions [13, 15, 16], and the limit (3) holds for a class of unbounded functions $F: \mathsf{X} \to \mathbb{R}$ [3, 49]. These conclusions are refined in [34] where in particular precise limit theory is obtained, generalizing the expansions of Bahadur and Ranga Rao for the partial sums of independent random variables [2, 4, 11]. Similar results are obtained in [33] for bounded functions on $\mathsf{X}$ under geometric ergodicity alone. Explicit, finite-time bounds have been obtained for uniformly ergodic chains in [21, 32].

Although most of the theory is based on assumptions on the Markov chain that are far stronger than geometric ergodicity, these conditions can be relaxed to obtain a weaker "pinned LDP" [41, 42]. Lower bounds can be obtained under essentially irreducibility alone [8, 29, 48].

The function $I: \mathbb{R} \to \mathbb{R}_+ \cup \{\infty\}$ appearing in (3) has many possible representations. In the limit theory of [3, 34, 41, 42] and the bounds obtained in [8, 29], the rate function is expressed as the convex dual

$$(4) \qquad I(c) = \sup_{a \in \mathbb{R}}[ca - \Lambda(a)], \qquad c \in \mathbb{R},$$

where the "pinned" log moment generating function is defined as

$$(5) \qquad \Lambda(a) = \lim_{n \to \infty} n^{-1} \log \mathsf{E}_\nu[\exp(naL_n(F))\mathbb{I}\{\Phi(n) \in C\}], \qquad a \in \mathbb{R},$$

with $C \subset \mathsf{X}$ a "small set" and $\nu$ a "small measure" (see discussion in Section 2.1). Under the assumptions imposed in this aforementioned work, the limit (5) exists, though it may be infinite, and is independent of the particular pair $(C, \nu)$ chosen (see [43], and the review in Section 2.1).

Sufficient as well as necessary conditions for the central limit theorem for Markov chains are presented in [20, 22, 23, 38, 43]. Much of the theory is based upon the *fundamental kernel*. Recall that a real-valued kernel $\widehat{P}$ on $\mathsf{X} \times \mathcal{B}$ is viewed as a linear operator, acting on functions $h: \mathsf{X} \to \mathbb{R}$ and probability measures $\mu$ on $\mathcal{B}$, via

$$(6) \qquad \widehat{P}h(\cdot) = \int_\mathsf{X} \widehat{P}(\cdot, dy)h(y) \quad \text{and} \quad \mu\widehat{P}(\cdot) = \int_\mathsf{X} \mu(dx)\widehat{P}(x, \cdot).$$

Under appropriate assumptions on $\mathbf{\Phi}$, the fundamental kernel can be expressed for an appropriate class of functions $F: \mathsf{X} \to \mathbb{R}$ via

$$ZF = \sum_{k=0}^{\infty} (P^k F - \pi(F)). \tag{7}$$

The following bilinear and quadratic forms are defined for measurable functions $F, G: \mathsf{X} \to \mathbb{R}$,

$$\langle\!\langle F, G \rangle\!\rangle := P(FG) - (PF)(PG), \qquad \mathcal{Q}(F) := P(F^2) - (PF)^2.$$

Under appropriate conditions, the asymptotic variance given in (2) can be expressed $\sigma^2(F) = \pi(\mathcal{Q}(ZF))$ (see [38], Theorem 17.5.3, and Proposition 2.1 below).

The purpose of the control-variate method is to reduce the variance of the standard estimator defined by

$$\phi_n := L_n(F), \qquad n \geq 1. \tag{8}$$

Suppose that there is a $\pi$-integrable function $H: \mathsf{X} \to \mathbb{R}$ with known mean. By normalization we can assume that $\pi(H) = 0$, and $L_n(F_\theta)$ is an asymptotically unbiased estimator of $\phi$ for each $\theta \in \mathbb{R}$ with $F_\theta := F - \theta H$.

The asymptotic variance of the controlled estimator is given by

$$\sigma^2(F_\theta) = \pi(\mathcal{Q}(ZF_\theta)) = \pi(\langle\!\langle ZF, ZF \rangle\!\rangle - 2\theta \langle\!\langle ZF, ZH \rangle\!\rangle + \theta^2 \langle\!\langle ZH, ZH \rangle\!\rangle).$$

Minimizing over $\theta \in \mathbb{R}$ gives the estimator with minimal asymptotic variance,

$$\theta^* = \frac{\pi(\langle\!\langle ZF, ZH \rangle\!\rangle)}{\pi(\langle\!\langle ZH, ZH \rangle\!\rangle)}.$$

See [17, 18, 36, 40] for more details and background on the general control-variate method.

An approach considered in [24, 25] is to consider functions of the form $H = J - PJ$, and choose $J$ so that it approximates the solution $\widehat{F}$ to *Poisson's equation*,

$$P\widehat{F} = \widehat{F} - F + \phi. \tag{9}$$

The idea is that if $J = \widehat{F}$, then the resulting controlled estimator with $\theta = 1$ has *zero* asymptotic variance.

This approach has been successfully applied in queueing models by taking $J$ equal to an associated *fluid value function*. The approach is provably effective in simple models [25], and numerical examples show dramatic variance reduction for more complex networks [26, 27]. Some theory to help explain the results of [26] is developed in [37] based on large deviation limit results contained in [33].

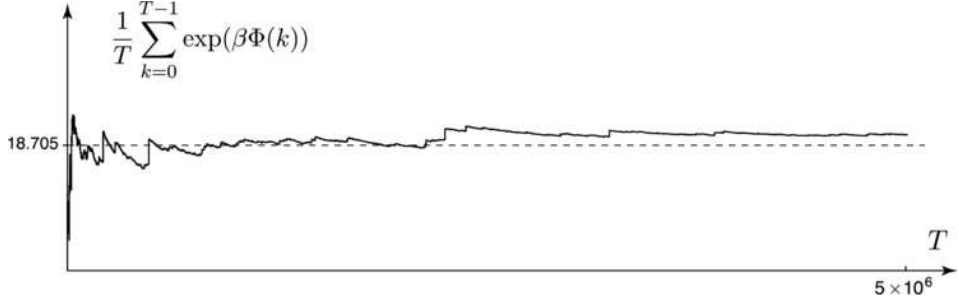

FIG. 1. *Monte Carlo estimates of $\phi := \pi(F)$ with $F(x) = e^{0.1x}$ for $x \in \mathbb{Z}_+$. The stochastic process $\mathbf{\Phi}$ is an $M/M/1$ queue initialized at zero, with load $\rho = 0.9$. After a transient period, the estimates are consistently larger than the steady-state mean of $\phi = (1 - \rho e^\beta)^{-1}(1 - \rho)$.*

For the sake of illustration consider the reflected random walk on $\mathbb{R}_+$, defined by the recursion

(10) $$\Phi(k+1) = [\Phi(k) + D(k+1)]_+, \qquad k \geq 0,$$

with $[x]_+ = \max(x, 0)$ for $x \in \mathbb{R}$, and $\mathbf{D}$ i.i.d. Consider first the special case in which $\mathbf{D}$ has common marginal distribution,

$$D(k) = \begin{cases} 1, & \text{with probability } \alpha, \\ -1, & \text{with probability } 1 - \alpha. \end{cases}$$

In this case $\mathbf{\Phi}$ is a discrete-time model of the $M/M/1$ queue, and the state space is then restricted to $\mathsf{X} = \mathbb{Z}_+$. It is assumed that $\alpha \in (0, \frac{1}{2})$ so that $\mathbf{\Phi}$ is a positive recurrent Markov chain on $\mathsf{X}$.

The invariant distribution $\pi$ is geometric, so there is little motivation to simulate. However, ignoring this issue momentarily, suppose we wish to estimate using simulation the steady-state mean of $F(x) = e^{\beta x}$ for a given $\beta > 0$.

Shown in Figure 1 are Monte Carlo estimates of the steady-state mean,

$$\phi := \sum \pi(i) F(i) = (1 - \rho) \sum \rho^i e^{\beta i},$$

where $\rho := \alpha/(1 - \alpha)$. In this simulation $\beta = 0.1$ and $\alpha = 9/19$, so that $\rho = 9/10$ and $\phi = (1 - \rho e^\beta)^{-1}(1 - \rho) \approx 18.705$. The Markov chain $\mathbf{\Phi}$ was initialized at zero, $\Phi(0) = 0$. The runlength in this simulation extended to $T = 5 \times 10^6$, yet the estimates are significantly larger than the steady-state mean over much of the run. The following proposition provides some explanation. A proof is provided in Section 3.2.

The existence of a nontrivial LDP depends upon structure of the sublevel sets of $F$, defined by $C_F(r) := \{x : F(x) \leq r\}$ for $r \geq 1$. This structure holds in Proposition 1.1(ii) since the sublevel sets are finite for each $r$.

PROPOSITION 1.1. *Consider the $M/M/1$ queue with $\rho = \alpha/(1-\alpha) < 1$.*

(i) *The Markov chain is geometrically ergodic, and its marginal distribution is geometric with parameter $\rho$.*

(ii) *Consider the function $F(x) = x$ or $F(x) = e^{\beta x}$, $x \in \mathbb{Z}_+$, for some fixed $\beta \in (0, |\log(\rho)|)$. The Monte Carlo estimates of the steady-state mean $\phi := \pi(F)$ are consistent, and there exists $\overline{c} < \phi$ such that the LDP (3) holds for any open set $O \subset (\overline{c}, \infty)$, and each initial condition $\Phi(0) = x \in \mathbb{Z}_+$. The convex rate function $I : [\overline{c}, \infty) \to \mathbb{R}_+$ is strictly positive on $[\overline{c}, \phi)$ and can be expressed as the convex dual (4). The rate function is identically zero on $[\phi, \infty)$. Consequently, we have for each initial condition $x \in \mathbb{Z}_+$,*

$$\lim_{n \to \infty} \frac{1}{n} \log(\mathsf{P}_x\{L_n(F) \leq c\}) = -I(c) < 0, \qquad c \in (\overline{c}, \phi),$$

*and*

$$\lim_{n \to \infty} \frac{1}{n} \log(\mathsf{P}_x\{L_n(F) \geq c\}) = 0, \qquad c \geq \phi.$$

In this example the chain is geometrically ergodic, so an LDP bound might not be surprising since an exact LDP holds when $F$ is bounded [33]. Section 3.2 contains a similar example in which analogous conclusions hold, yet $\mathbf{\Phi}$ is not geometrically ergodic, and $\{\phi_n\}$ does not even satisfy the CLT. Moreover, in this example a control-variate is constructed to obtain a *pair* of estimators giving upper *and* lower confidence bounds.

The remainder of the paper is organized as follows. Section 2 contains a statement and proof of the most important conclusion in this paper, Theorem 2.2, which establishes the LDP for a general class of functions on X. Section 2.1 contains a survey of spectral theory for Markov chains, following [3, 33, 34, 39]. A new criterion for the existence of a spectral gap is presented in Section 2.2, which is the main ingredient in the proof of Theorem 2.2.

Applications of Theorem 2.2 to the construction and analysis of control-variates are contained in Section 3.1. The simulation algorithm proposed in Section 3.1 is shown to satisfy exact upper and lower LDP bounds. This result is illustrated in Section 3.2 using the reflected random walk (10). Conclusions are contained in Section 4.

**2. One-sided large deviation asymptotics.** Throughout the paper it is assumed that $\mathbf{\Phi}$ is positive Harris recurrent and aperiodic. Equivalently, there is a unique invariant probability distribution $\pi$ on $\mathcal{B}$ such that, for any $A \in \mathcal{B}$ satisfying $\pi(A) > 0$, and any initial condition $x$,

$$\lim_{k \to \infty} \|P^k(x, \cdot) - \pi(\cdot)\| = 0, \qquad x \in \mathsf{X},$$

where $\|\cdot\|$ denotes the total-variation norm ([38], Theorem 13.0.1). We denote by $\mathcal{B}^+$ the set of $A \in \mathcal{B}$ satisfying $\pi(A) > 0$. We write $f \in \mathcal{B}^+$ if $f : \mathsf{X} \to \mathbb{R}_+$ is a measurable function with $\int f \, d\pi > 0$.

A measurable function $s : \mathsf{X} \to \mathbb{R}_+$ and a probability measure $\nu$ on $\mathcal{B}$ are called *small* if for some $n \geq 1$ we have

$$(11) \qquad P^n(x, A) \geq s(x)\nu(A), \qquad x \in \mathsf{X}, A \in \mathcal{B}.$$

A set $C$ is called *small* if $s = \epsilon \mathbb{I}_C$ is a small function for some positive $\epsilon$.

The following *Lyapunov drift condition* is assumed throughout the paper. Given any measurable function $F : \mathsf{X} \to \mathbb{R}$ satisfying $\pi(|F|) < \infty$, we can construct a solution to (V3) with $W = 1 + |F|$ by applying [38], Theorem 14.2.3. The set $C_V$ on which $V$ is finite is absorbing, so that the chain can be restricted to this set; see [38], Proposition 4.2.3.

$$(\text{V3}) \qquad \begin{cases} \text{For a function } W : \mathsf{X} \to [1, \infty), \\ \text{a small set } C \subset \mathsf{X}, \text{ and a constant } b < \infty, \\ PV \leq V - W + b\mathbb{I}_C \qquad \text{on } C_V := \{x : V(x) < \infty\}. \end{cases}$$

For a given function $W : \mathsf{X} \to [1, \infty)$ the weighted $L_\infty$-norm is defined as

$$\|h\|_W := \sup_x |h(x)|/W(x),$$

and $L_\infty^W$ denotes the set of all measurable functions $h : \mathsf{X} \to \mathbb{R}$ for which this norm is finite (see [28, 30, 31, 33, 34, 38]). The supremum norm $\|\cdot\|_\infty$ is precisely $\|\cdot\|_W$ with $W \equiv 1$. Two functions $W, W' : \mathsf{X} \to [1, \infty)$ are called *equivalent* if they generate the same function space, that is,

$$W' \in L_\infty^W \quad \text{and} \quad W \in L_\infty^{W'}.$$

The set of finite measures on $\mathcal{B}$ is denoted $\mathcal{M}$; the set $\mathcal{M}_1 \subset \mathcal{M}$ denotes probability measures on $\mathcal{B}$; $\mathcal{M}^W \subset \mathcal{M}$ denotes measures satisfying $\mu(W) < \infty$; and $\mathcal{M}_1^W = \mathcal{M}^W \cap \mathcal{M}_1$.

The convergence results in parts (i) and (ii) of Proposition 2.1 are contained in the $f$-norm ergodic theorem (Theorem 14.0.1) of [38]. The interpretation of the sum as a version of the fundamental kernel is contained in [38], Theorem 17.4.2.

Part (iii) follows from [38], Theorem 17.4.4, and (iv) is contained in [38], Theorem 16.0.1.

PROPOSITION 2.1. *Suppose that $\mathbf{\Phi}$ is $\psi$-irreducible and aperiodic and* (V3) *holds with $V$ everywhere finite. Then:*

(i) *The chain is positive Harris recurrent with $\pi(W) < \infty$, and we have*

$$\lim_{k \to \infty} \|P^k(x, \cdot) - \pi(\cdot)\|_W = 0, \qquad x \in \mathsf{X}.$$

*Moreover, the fundamental kernel $Z : L_\infty^W \to L_\infty^V$ exists as a bounded linear operator.*

(ii) *If $\pi(V) < \infty$, then the chain $\Phi$ is called $W$-regular of degree* 2. *In this case the fundamental kernel can be expressed as the sum* (7). *The sum converges in the induced operator norm from $L_\infty^W$ to $L_\infty^V$.*

(iii) *If $\pi(WV) < \infty$, then the CLT* (2) *holds for each $F \in L_\infty^W$.*

(iv) *If $W \geq \varepsilon_0 V$ for some $\varepsilon_0 > 0$ then we say that* (V4) *holds. In this case $\Phi$ is $V$-uniformly ergodic: for some constants $b_0 < \infty$, $r_0 > 1$,*

$$\sum_{k=0}^\infty r_0^k \|P^k(x,\cdot) - \pi(\cdot)\|_V \leq b_0(V(x) + 1), \qquad x \in \mathsf{X}.$$

We list below some other definitions for a given measurable function $F: \mathsf{X} \to \mathbb{R}$: The function is called

*Degenerate* if there is a measurable function $H: \mathsf{X} \to \mathbb{R}$ such that when $\Phi(0) \sim \pi$, $F(\Phi(k+1)) - F(\Phi(k)) = H(\Phi(k))$ a.s. for $k \geq 0$. Under appropriate bounds on $H$ this implies that the asymptotic variance of $\phi_n$ is equal to zero. A converse is provided in [34], Lemma 4.12, based on [33], Proposition 2.4.

*Lattice* if there are $h > 0$ and $0 \leq d < h$, such that

(12) $$\frac{F(x) - d}{h} \text{ is an integer}, \qquad x \in \mathsf{X}.$$

If there exists a lattice function $F_\ell$ such that $F - F_\ell$ is degenerate, then $F$ is called *almost-lattice*. Otherwise, $F$ is called *strongly nonlattice*.

*Near-monotone* if $\inf_{x \in \mathsf{X}} F(x) > -\infty$, and the sublevel set $C_F(r) := \{x \in \mathsf{X} : F(x) \leq r\}$ is small or empty for each $r < \|F_+\|_\infty$, where $F_+ = \max(F, 0)$ and $\|\cdot\|_\infty$ denotes the supremum norm.

Large deviation bounds are obtained in [3] for countable state-space chains under the assumption that $F$ is near-monotone, and $\Lambda(F) < \|F\|_\infty$, where $\Lambda(F)$ is defined in (5) using $a = 1$. These assumptions are far stronger than geometric ergodicity when $F$ is unbounded. The results of [3] are strengthened and generalized to general state-space chains and processes in [33, 34].

By restricting the range of $c$ in (3) we can relax the geometric ergodicity assumption. The proof of Theorem 2.2 is included at the end of this section.

The most important assumption in Theorem 2.2 is the constraint (13). To interpret this condition, consider first the countable state-space case. If $V$ and $F$ have finite sublevel sets and $V$ is unbounded, then this condition is immediate since $C_F(r)$ and $C_V(r_0)$ are each finite for $r < \|F\|_\infty$ and $r_0 < \infty$, and $C_V(r_0) \uparrow \mathsf{X}$ as $r_0 \uparrow \infty$.

For general state-space models the set $C_V(r_0)$ is always $F$-*regular* for any finite $r_0$, and hence small, by Theorem 14.0.1 combined with Proposition 14.1.2 of [38]. In this way we can interpret (13) as simultaneously a relaxation and strengthening of the near-monotone condition.

THEOREM 2.2. *Suppose that* (V3) *holds with $V$ everywhere finite, and that $F \in L_\infty^W$ is a nondegenerate function satisfying $\pi(F) = 0$. Suppose moreover that the sublevel set $C_F(r)$ satisfies for some $r > \phi = 0$, $r_0 < \infty$,*

$$(13) \qquad C_F(r) \subset C_V(r_0).$$

*Then, there exists $\overline{c}_0 < \phi$ and a smooth convex function $I : (\overline{c}_0, \phi) \to (0, \infty)$ such that:*

(i) *The LDP* (3) *holds for each initial condition $x \in \mathsf{X}$ and each $c \in (\overline{c}_0, \phi)$.*

(ii) *If $F$ is strongly nonlattice, then for each $c \in (\overline{c}_0, \phi)$, there exists a bounded function $g_c : \mathsf{X} \to (0, \infty)$, such that for each initial condition $x \in \mathsf{X}$,*

$$(14) \qquad \mathsf{P}_x\{L_n(F) \leq c\} \sim \frac{g_c(x)}{\sqrt{2\pi n}} e^{-nI(c)}, \qquad n \to \infty.$$

The LDP asymptotics described in Theorem 2.2 are based on the spectral theory of a positive semigroup obtained from the function $F$ to be simulated. The definitions presented in Section 2.1 are taken from [3, 33, 34, 38, 39, 43].

2.1. *Positive semigroups.* Consider now a positive kernel $\widehat{P}$ on $\mathsf{X} \times \mathcal{B}$. It is assumed that the semigroup $\{\widehat{P}^k : k \geq 0\}$ is $\psi$-*irreducible*,

$$\sum_{k=0}^{\infty} \widehat{P}^k(x, A) > 0, \qquad x \in \mathsf{X}, A \in \mathcal{B}^+,$$

and also *aperiodic*,

$$\liminf_{k \to \infty} \mathbb{I}\{\widehat{P}^k(x, A) > 0\} = 1 \qquad \text{for each } x \in \mathsf{X}, A \in \mathcal{B}^+.$$

For a $\psi$-irreducible kernel there exists a function $s \in \mathcal{B}^+$, a probability measure $\nu$ on $\mathcal{B}$, and $n_0 \geq 1$ such that $\widehat{P}^{n_0} \geq s \otimes \nu$. The function $s$ and measure $\nu$ are called $\widehat{P}$-*small*, generalizing the definition for a probabilistic kernel $P$.

Based on a given function $h : \mathsf{X} \to (0, \infty)$ we consider in Section 2.2 the two positive kernels,
the *scaled kernel*: $P_h := I_h P$, or equivalently,

$$P_h(x, A) := h(x) P(x, A), \qquad x \in \mathsf{X}, A \in \mathcal{B},$$

the *twisted kernel*: $\check{P}_h := I_{Ph}^{-1} P I_h$, or equivalently,

$$(15) \qquad \check{P}_h(x, A) = \frac{\int_A P(x, dy) h(y)}{Ph(x)}, \qquad x \in \mathsf{X}, A \in \mathcal{B}.$$

The twisted kernel is *probabilistic*, so that $\check{P}_h(x, \mathsf{X}) = 1$ for all $x$, provided $Ph(x) < \infty$ for all $x \in \mathsf{X}$.

For any $\psi$-irreducible and aperiodic semigroup, the *generalized principal eigenvalue* (g.p.e.) is defined as $\lambda = e^\Lambda$, where $\Lambda \in [-\infty, \infty]$ is the limit,

(16) $$\Lambda := \lim_{n \to \infty} n^{-1} \log(\nu \widehat{P}^n s).$$

The limit is independent of the particular small function $s \in \mathcal{B}^+$ and small measure $\nu$ chosen. If $\Lambda$ is finite, then there is an associated eigenfunction $h : \mathsf{X} \to (0, \infty]$ satisfying $h(x) < \infty$ a.e. $[\psi]$, and

$$\widehat{P} h \leq \lambda h.$$

This is an equality provided $\widehat{P}$ is $\lambda$-*recurrent*,

$$\sum_{k=0}^\infty \lambda^{-k} \nu \widehat{P}^k s = \infty.$$

See [9, 33, 43] for further discussion.

For a given weighting function $v : \mathsf{X} \to [1, \infty)$, the induced operator norm of $\widehat{P}$ is

(17) $$\|\|\widehat{P}\|\|_v := \sup\left\{ \frac{\|\widehat{P} h\|_v}{\|h\|_v} : h \in L_\infty^v, \|h\|_v \neq 0 \right\}.$$

The *spectrum* $\mathcal{S}(\widehat{P}) \subset \mathbb{C}$ of $\widehat{P}$ is the set of $z \in \mathbb{C}$ such that the inverse $[Iz - \widehat{P}]^{-1}$ does not exist as a bounded linear operator on $L_\infty^v$.

The *spectral radius* of the semigroup $\{\widehat{P}^t\}$ is expressed $\xi = e^\Xi$, where

(18) $$\Xi := \lim_{k \to \infty} k^{-1} \log(\|\|\widehat{P}^k\|\|_v).$$

We say that $\widehat{P}$ is $v$-*uniform* if the spectral radius $\xi$ is finite, and there exists $h \in L_\infty^v$, $\mu \in \mathcal{M}_1^v$, such that

$$\sup_{|z| \geq \xi} \|\|[Iz - (\widehat{P} - h \otimes \mu)]^{-1}\|\|_v < \infty.$$

When $v \equiv 1$ we drop the qualification and simply say that $\widehat{P}$ is *uniform*.

If the kernel is $v$-uniform, then it admits a spectral gap, and the generalized principal eigenvalue coincides with $\xi$. Moreover, the eigenfunction $h$ satisfies $h \in L_\infty^v$, and the eigenfunction equation $\widehat{P} h = \lambda h$ holds ([34], Proposition 2.9).

2.2. *Multiplicative mean-ergodic theorem.* The multiplicative mean-ergodic theorem contained in Theorem 2.3 is the basis of LDP asymptotics for the partial sums [3, 33, 34].

The proof of Theorem 2.3 is identical to the proof of Theorem 3.4 in [34]. The idea of the proof of (20) is as follows: The twisted kernel $\check{P}_h$ is $\check{v}$-uniform, where $h = \check{f}$ is an eigenfunction and $\check{v} := v/h$, since the twisted kernel is simply a scaling and similarity transformation of $P_f$ with $f = e^F$. Since each of the twisted kernels is probabilistic, this implies that $\check{P}_h$ is the transition kernel for a $\check{v}$-uniformly ergodic Markov chain (see [33], Corollary 4.7 or [34], Proposition 2.11 for finer results). The multiplicative mean-ergodic theorem is a consequence of this mean-ergodic theorem for the twisted chain, and the representation

$$P_f^n(x, A) = \mathsf{E}_x[\exp(nL_n(F))\mathbb{I}\{\Phi(n) \in A\}], \qquad x \in \mathsf{X}, A \in \mathcal{B}.$$

An explicit formula for the eigenfunction $\check{f}$ is given in (28).

THEOREM 2.3. *Suppose that $F: \mathsf{X} \to \mathbb{R}$ is measurable; its g.p.e. $\lambda$ is finite; and that $P_f$ is $v$-uniform, with $f = e^F$. Then, there exists a measure $\check{\mu} \in \mathcal{M}_1^v$ and a function $\check{f} \in L_\infty^v$ satisfying the eigenfunction equations $\widehat{P}\check{f} = \lambda\check{f}$, $\check{\mu}\widehat{P} = \lambda\check{\mu}$, with normalization,*

$$\check{\mu}(\check{f}) = \check{\mu}(\mathsf{X}) = 1, \tag{19}$$

*and these are the unique solutions. Moreover, the following multiplicative mean-ergodic theorem holds: For some $b_0 > 0$, $b_1 < \infty$ and all $x \in \mathsf{X}$, $n \geq 1$,*

$$\left| \mathsf{E}_x\left[ \exp\left( \sum_{k=0}^{n-1} F(\Phi(k)) - n\Lambda \right) \right] - \check{f}(x) \right| \leq b_1 e^{-b_0 n} v(x). \tag{20}$$

In the series of results that follow we present sufficient conditions for uniformity of a scaled kernel. It will be convenient to consider a family of kernels $\{P_a : a \in (0, 1]\}$ where $P_a = P_{f_a}$ is the scaled kernel defined with $f_a = e^{aF}$, and $F: \mathsf{X} \to \mathbb{R}$ a given measurable function. A family of resolvent kernels is defined by

$$R_a := \sum_{k=0}^\infty 2^{-k-1} P_a^k, \qquad a \in (0, 1], \tag{21}$$

and $R$ denotes the kernel obtained when $a = 0$ so that $f_a \equiv 1$. It can be shown that $v$-uniformity of $R_a$ is equivalent to $v$-uniformity of $P_a$ when $P$ is aperiodic and $\lambda_a < 2$, and in this case $\lambda_a$ is the g.p.e. for $P_a$ if and only if $\gamma_a = (2 - \lambda_a)^{-1}$ is the g.p.e. for $R_a$.

We assume throughout that the function $F$ is normalized so that $\pi(F) = 0$. It then follows from the definition (16), Jensen's inequality and the mean-ergodic theorem that $\lambda_a, \gamma_a \in [1, \infty]$ for each $a \in \mathbb{R}$.

Under the assumptions of Theorem 2.4 the kernel $R_a$ is a bounded linear operator on $L_\infty^{v_a}$ for $a \in [0, 1]$, where $v_a = e^{aV}$. The first bound in (22) is analogous to Condition (V4) of [38]. These two bounds are equivalent to geometric ergodicity for the Markov chain in the special case $F \equiv 0$. For general $F$ this is not true, as we shall see in Theorem 2.6.

THEOREM 2.4. *Suppose that $F:\mathsf{X} \to \mathbb{R}$ is a given measurable function satisfying $\pi(|F|) < \infty$ and $\pi(F) = 0$. Suppose that there exist $\bar{\lambda}_1 < 1$, a function $v:\mathsf{X} \to [1,\infty]$ that is not everywhere infinite, a function $s:\mathsf{X} \to \mathbb{R}_+$, a probability distribution $\nu$ on $\mathcal{B}$ and a constant $b < \infty$ satisfying the bounds*

$$P_f v \leq \bar{\lambda}_1 v + bs \quad \text{and}$$
(22)
$$R_a \geq s \otimes \nu \qquad \text{for all } 0 \leq a \leq 1.$$

*Then, $v(x) < \infty$ a.e. $[\pi]$, and there exists $\bar{a} \in (0,1)$ such that $P_a$ is $v_a$-uniform for all $a \in (0, \bar{a})$.*

PROOF. To show that $R_a$ is $v_a$-uniform we prove that $\|\!|G_a|\!\|_{v_a} < \infty$, where $G_a$ denotes the potential kernel,

(23)
$$G_a = \sum_{k=0}^{\infty} \gamma_a^{-k-1}[R_a - s \otimes \nu]^k,$$

and $\gamma_a = (2 - \lambda_a)^{-1}$ is the g.p.e. for $R_a$.

Jensen's inequality implies the following family of bounds:

$$P_a v_a \leq (\bar{\lambda}_1 v + bs)^a \leq \bar{\lambda}_a v_a + a\bar{\lambda}_{a-1} bs, \qquad 0 \leq a \leq 1,$$

where $\bar{\lambda}_t := (\bar{\lambda}_1)^t$ for any $t \in \mathbb{R}$. Moreover, the resolvent equation holds: $P_a R_a = R_a P_a = 2R_a - I$. This combined with the bound on $P_a v_a$ gives

(24)
$$2 R_a v_a - v_a = R_a P_a v_a \leq R_a[\bar{\lambda}_a v_a + a\bar{\lambda}_{a-1} bs],$$

which on rearranging terms implies the bound

(25)
$$R_a v_a \leq \bar{\gamma}_a v_a + a b_a R_a s,$$

with $\bar{\gamma}_a = (2 - \bar{\lambda}_a)^{-1}$, and $b_a = b\bar{\lambda}_{a-1}\bar{\gamma}_a$. We evidently have $\bar{\gamma}_a < 1$ and $\bar{\lambda}_a < 1$, so that $b_a \leq b/\bar{\lambda}_1$ for $a \in [0,1]$.

Define $v_a' = (1 + ab_2)v_a - ab_2 s$ where $b_2 > b/\bar{\lambda}_1$ is fixed. This function is equivalent to $v_a$, and from (25),

$$R_a v_a' \leq (1 + ab_2)[\bar{\gamma}_a v_a + ab_a R_a s] - ab_2 R_a s.$$

For $a > 0$ sufficiently small we have $b_2 \geq (1 + ab_2)b/\bar{\lambda}_1 \geq (1 + ab_2)b_a$. Consequently, for such $a$,

$$R_a v_a' \leq (1 + ab_2)\bar{\gamma}_a v_a = \bar{\gamma}_a v_a' + ab_2 \bar{\gamma}_a s,$$

and on subtracting the function $\nu(v_a')s$ from each side this gives

$$[R_a - s \otimes \nu]v_a' \leq \bar{\gamma}_a v_a' - (\nu(v_a') - ab_2\bar{\gamma}_a)s.$$

Decreasing $a$ still further we can assume that $(\nu(v_a') - ab_2\bar{\gamma}_a) \geq 0$. We conclude that there exists $\bar{a} > 0$ such that with $\delta_a := 1 - \bar{\gamma}_a$,

$$[R_a - s \otimes \nu]v_a' \leq v_a' - \delta_a v_a', \qquad a \in (0, \bar{a}].$$

Iterating this bound gives

$$[R_a - s \otimes \nu]^n v'_a \leq v'_a - \delta_a \sum_{k=0}^{n-1}[R_a - s \otimes \nu]^k v'_a, \qquad n \geq 1,$$

and hence, $\sum_{k=0}^{\infty}[R_a - s \otimes \nu]^k v'_a \leq v'_a$, which implies the final bound,

$$\|\|[Iz - (R_a - s \otimes \nu)]^{-1}\|\|_{v'_a} \leq \delta_a^{-1}, \qquad a \in (0, \bar{a}], \ |z| \geq 1.$$

This completes the proof that $\|\|G_a\|\|_{v_a} < \infty$ since $\gamma_a \geq 1$, and $v_a$ is equivalent to $v'_a$ for $a \in (0, \bar{a}]$. □

The following result provides a simple criterion that guarantees the existence of $s, \nu$ satisfying the minorization condition in (22).

PROPOSITION 2.5. *Suppose that* (V3) *holds, and that* $F \in L_\infty^W$. *Then, for each* $r \geq 1$, $a \in \mathbb{R}$, *the set* $C = C_V(r) := \{x : V(x) \leq r\}$ *is small for the positive kernel* $P_a$. *Moreover, we have the following uniform bounds: For each* $a_0 > 0$, $r \geq 1$, *there exist* $\varepsilon_0 > 0$, $n_0 \geq 1$ *and a probability distribution* $\nu_0 \in \mathcal{M}_1$ *such that,*

$$P_a^{n_0}(x, A) \geq \varepsilon_0 \nu_0(A), \qquad x \in C, \ a \in [-a_0, a_0].$$

PROOF. Fix $a_0 > 0$, $r \geq 1$, set $F_0(x) = a_0\|F\|_W W(x)$, $x \in \mathsf{X}$, and define $\widehat{P} := I_{e^{-F_0}} P$. This is simply the scaled kernel $P_h$ with $h = e^{-F_0}$. A minorization condition obtained for the kernel $\widehat{P}$ will imply the desired uniform bounds since $P_a \geq \widehat{P}$ for $a \in [-a_0, a_0]$.

For any $A \in \mathcal{B}$ and $x \in \mathsf{X}$ we have by Jensen's inequality,

$$\widehat{P}^n(x, A) = \left\{ P^n(x, A)^{-1} \mathsf{E}_x\left[\exp\left(-\sum_{k=0}^{n-1} F_0(\Phi(k))\right) \mathbb{I}\{\Phi(n) \in A\}\right] \right\} P^n(x, A)$$
(26)
$$\geq \exp\left\{ P^n(x, A)^{-1} \mathsf{E}_x\left[\left(-\sum_{k=0}^{n-1} F_0(\Phi(k))\right) \mathbb{I}\{\Phi(n) \in A\}\right] \right\} P^n(x, A).$$

Under (V3) the following bound holds: $0 \leq P^n V \leq V + nb - \sum_{k=0}^{n-1} P^k W$, so that

$$\mathsf{E}_x\left[\sum_{k=0}^{n-1} F_0(\Phi(k))\right] \leq a_0\|F\|_W[V(x) + nb], \qquad x \in \mathsf{X}, n \geq 1.$$

Consequently, from (26),

(27)
$$\widehat{P}^n(x, A) \geq \exp\{-P^n(x, A)^{-1} a_0\|F\|_W[V(x) + nb]\} P^n(x, A),$$
$$x \in \mathsf{X}, n \geq 1.$$

This shows that the semigroup $\{\widehat{P}^n : n \geq 1\}$ is $\pi$-irreducible, in the sense that $\sum_n \widehat{P}^n(x,A) > 0$ for each $x$ whenever $\pi(A) > 0$. Let $A \in \mathcal{B}^+$ be any fixed small set for $\widehat{P}$; there is a probability distribution $\nu_0$, $\varepsilon > 0$ and an integer $m \geq 1$ such that

$$\widehat{P}^m(x,B) \geq \varepsilon \nu_0(B), \qquad x \in A, \ B \in \mathcal{B}.$$

Choose $n \geq 1$, $\delta > 0$ such that $P^n(x,A) \geq \delta$ for $x \in C$. This is possible since the set $C$ is $W$-regular, and hence small ([38], Theorem 14.2.3). It follows from (27) that

$$\widehat{P}^n(x,A) \geq \delta_r := \exp\{-\delta^{-1} a_0 \|F\|_W [r+nb]\} \delta, \qquad x \in C,$$

and hence,

$$\widehat{P}^{n+m}(x,B) \geq \delta_r \varepsilon \nu_0(B), \qquad x \in C, B \in \mathcal{B}.$$

This completes the proof with $n_0 = m + n$ and $\varepsilon_0 = \delta_r \varepsilon$.  □

Up to now it appears that uniformity is a tremendously strong assumption on the scaled kernel $P_f$ since the implications of uniformity are so strong. However, under the assumptions of Theorem 2.2 we can establish uniformity of $P_a$ for a range of *nonpositive* $a$, even though $\boldsymbol{\Phi}$ is only positive Harris recurrent.

Recall that $\widehat{P}$ is called *uniform* if it is $v$-uniform with $v \equiv 1$.

THEOREM 2.6.  *Suppose that* (V3) *holds with $V$ everywhere finite, and suppose that the function $F \in L_\infty^W$ satisfies* (13) *for some $r > 0$ and $r_0 < \infty$, with $\phi = \pi(F) = 0$. Then, there exists $\bar{a} < 0$ such that:*

  (i) $P_a$ *is uniform for each* $a \in (\bar{a}, 0)$.
  (ii) *The eigenfunctions* $\{\check{f}_a : a \in (\bar{a}, 0)\} \subset L_\infty$, *normalized so that* $\nu(\check{f}_a) = 1$ *for some small measure $\nu$ and each $a$, are uniformly bounded:*

$$\sup_{\substack{\bar{a} < a < 0 \\ x \in \mathsf{X}}} \check{f}_a(x) < \infty.$$

  (iii) *Define* $\check{f}_a' := \frac{d}{da} \check{f}_a$ *for* $a \in (\bar{a}, 0)$, *with* $\{\check{f}_a\}$ *normalized as in* (ii). *These functions are uniformly bounded in norm:*

$$\sup_{\bar{a} < a < 0} \|\check{f}_a'\|_V < \infty.$$

  (iv) $\Lambda$ *is convex and analytic on* $(\bar{a}, 0)$, *and* $\lim_{a \uparrow 0} \frac{d}{da} \Lambda(a) = 0$.

Define the twisted kernel by $\check{P}_a := \check{P}_{\check{f}_a}$, where $\check{f}_a$ is an eigenfunction that exists for $P_a$. We have noted prior to Theorem 2.3 that $\check{P}_a$ is the transition kernel for a $\check{v}_a$-uniformly ergodic Markov chain when $P_a$ is $v_a$-uniform. As in [33], Proposition 4.9 one can verify that each of the functions

$$\widehat{F}_a = \frac{d}{da}\log(\check{f}_a) = \check{f}'_a/\check{f}_a, \qquad a \in (\bar{a}, 0),$$

solves Poisson's equation for the corresponding twisted kernel,

$$\check{P}_a \widehat{F}_a = \widehat{F}_a - F + \phi(a),$$

where $\phi(a) = \frac{d}{da}\Lambda(a)$ is the steady-state mean of $F$ for the twisted kernel. Poisson's equation for $\check{P}_a$ is used to establish versions of Theorem 2.4(iv) in the papers [33, 34]. This technique cannot be applied here since we do not know if $\Lambda$ is bounded or smooth for positive $a$.

The proof of Theorem 2.2 is performed in the remainder of this subsection through a series of steps.

We see in Lemma 2.7 that part (i) of Theorem 2.2 follows quickly from Theorem 2.4. Recall that $G_a = [\gamma_a I - (R_a - s \otimes \nu)]^{-1}$ is the potential kernel previously defined in (23).

LEMMA 2.7. *If the assumptions of Theorem 2.6 hold, then there exists $\bar{a}_0 < 0$ such that $P_a$ is uniform for each $a \in (\bar{a}_0, 0)$. The unique eigenfunction $\check{f}_a \in L_\infty$ satisfying $\nu(\check{f}_a) = 1$ can be expressed*

(28) $$\check{f}_a = G_a s.$$

PROOF. Set $G = -F$, $g = e^G$ and $g_a = e^{aG}$ for $a \in \mathbb{R}$. Note that uniformity of $P_{g_a}$ is equivalent to uniformity of $P_{-a}$ for any $a$.

Define $v \equiv 1$, $\delta = r - \phi$ and $b = \exp(\sup_{x \in C_F(r)} |G(x)|) < \infty$, so that the following bound holds:

$$P_g v = e^{-F} \leq e^{-\delta} v + b \mathbb{I}_{C_F(r)}.$$

Moreover, Proposition 2.5 implies that the minorization condition in (22) holds with $F$ replaced by $G$ in the definition of $R_a$. Uniformity of $P_{g_a}$ for sufficiently small $a > 0$ thus follows from Theorem 2.4 and the fact that $v$ is bounded.

The representation (28) follows from Proposition 2.8 of [34] (see also [43]). □

The difficult part of the proof of Theorem 2.6 is to establish convergence of $\phi(a)$ to $\phi$ as $a \uparrow 0$. The proof is based on consideration of the following scaled kernel to bound $P_a$.

For a given small measure $\nu \in \mathcal{M}_1$, let $s = \varepsilon_v \mathbb{I}_{C_V(r_0)}$ with $\varepsilon_v > 0$ chosen so that $R \geq s \otimes \nu$. For a fixed $\varepsilon > 0$, define the scaled kernel,

$$\bar{P}(x, A) := \exp(\varepsilon_v \varepsilon \mathbb{I}_{C_V(r_0)}(x)) P(x, A), \qquad x \in \mathsf{X}, A \in \mathcal{B},$$

and the resolvent and potential kernels,

$$\bar{R} := \sum_{k=0}^{\infty} 2^{-k-1} \bar{P}^k, \qquad \bar{G} := \sum_{k=0}^{\infty} [\bar{R} - s \otimes \nu]^k.$$

LEMMA 2.8. *Suppose that the assumptions of Theorem 2.6 hold. Then, there exists $\varepsilon > 0$ such that:*

(i) $\|\bar{G}s\|_\infty < \infty$,
(ii) $\|\bar{G}\bar{R}s\|_\infty < \infty$,
(iii) $\|\bar{G}W\|_V < \infty$.

PROOF. To see (i) we write

$$\bar{P}\mathbf{1} = 1 + (e^{\varepsilon \varepsilon_v} - 1) \mathbb{I}_{C_V(r_0)} = 1 + \varepsilon_v^{-1} (e^{\varepsilon \varepsilon_v} - 1)s.$$

Let $h = 1 - \varepsilon_v^{-1}(e^{\varepsilon \varepsilon_v} - 1)s$, and choose $\varepsilon > 0$ so this is strictly positive everywhere. From the resolvent equation as in (24) we have

$$\bar{R}\mathbf{1} = 1 + \varepsilon_v^{-1}(e^{\varepsilon \varepsilon_v} - 1)\bar{R}s,$$

and hence

$$\bar{R}h = 1 = h + \varepsilon_v^{-1}(e^{\varepsilon \varepsilon_v} - 1)s.$$

On subtracting $[s \otimes \nu]h$ from both sides we then obtain

$$(\bar{R} - s \otimes \nu)h = h - \delta_h s,$$

where $\delta_h = \nu(h) - \varepsilon_v^{-1}(e^{\varepsilon \varepsilon_v} - 1)$. By reducing $\varepsilon > 0$ we can assume that $\delta_h > 0$.

Exactly as in the proof of Theorem 2.4 we conclude that

$$\bar{G}s := \sum_{k=0}^{\infty} [\bar{R} - s \otimes \nu]^k s \leq \delta_h^{-1} h,$$

which establishes the uniform bound in (i).

Part (ii) follows from (i) and the identity,

(29) $$\bar{G}\bar{R} = \bar{G}[\bar{R} - s \otimes \nu] + \bar{G}[s \otimes \nu] = \bar{G} - I + \bar{G}[s \otimes \nu],$$

so that $\bar{G}\bar{R}s \leq [\bar{G} + (\bar{G}s) \otimes \nu]s = (1 + \nu(s))\bar{G}s$.

To see (iii) we note that we can assume without loss of generality that the set $C$ in (V3) is equal to the set $C_V(r_0)$ used here by applying [38], Theorem 14.2.3. Under this transformation we obtain

$$\bar{P}V \leq V - W + (e^{\varepsilon s} - 1)V + (e^{\varepsilon \varepsilon_v} - 1)b \mathbb{I}_C$$
$$\leq V - W + b_v s,$$

where $b_v = \varepsilon_v^{-1}(e^{\varepsilon \varepsilon_v} - 1)(b + r_0)$.

From the resolvent equation again this implies the bound

$$\bar{R}V \leq V - \bar{R}W + b_v \bar{R}s,$$

and then through familar arguments, $\bar{G}\bar{R}W \leq V + b_v \bar{G}\bar{R}s$. This bound combined with the identity (29) completes the proof of (iii). □

With this value of $\varepsilon > 0$ fixed in the definition of $\bar{R}$, and given $\bar{a}_0 < 0$ from Lemma 2.7, we now identify the lower bound $\bar{a}$:

LEMMA 2.9. *Suppose that the assumptions of Theorem 2.6 hold. Then, there exists $\bar{a} \in (\bar{a}_0, 0)$ such that for any $\bar{a} < a \leq 0$,*

$$\lambda_a^{-1} P_a \leq P_a \leq \bar{P} \quad \text{and} \quad \check{f}_a \leq \bar{f} := \bar{G}s,$$

*where $\check{f}_a$ is defined in (28).*

PROOF. The bound $\lambda_a^{-1} P_a \leq P_a$ holds since $\lambda_a \geq 1$.

To see that $P_a \leq \bar{P}$, rewrite this bound as $f_a(x) \leq \exp(\varepsilon_v \varepsilon \mathbb{I}_{C_V(r_0)}(x))$, or on taking logarithms,

(30) $$aF(x) \leq \varepsilon_v \varepsilon \mathbb{I}_{C_V(r_0)}(x), \qquad x \in \mathsf{X}.$$

Letting $b_0 = \sup_{x \in C_0(F)} |F(x)|$ with $C_0(F) = \{x \in \mathsf{X} : F(x) \leq 0\}$, and applying the bound (13) gives,

$$aF(x) \leq |a| b_0 \mathbb{I}_{C_0(F)} \leq |a| b_0 \mathbb{I}_{C_V(r_0)}, \qquad a < 0.$$

This shows that (30) holds for $a \in [\bar{a}, 0)$ with

$$\bar{a} := \max(\bar{a}_0, -\varepsilon_v \varepsilon b_0^{-1}).$$

The second bound $\check{f}_a \leq \bar{f}$ follows immediately from the first since

$$\gamma_a^{-n-1}(R_a - s \otimes \nu)^n \leq (R_a - s \otimes \nu)^n \leq (\bar{R} - s \otimes \nu)^n, \qquad n \geq 0. \quad \square$$

The previous two lemmas lead to a proof of Theorem 2.6(iii):

LEMMA 2.10. *Under the assumptions of Theorem 2.6 the function $\check{f}'_a$ can be expressed*

$$\check{f}'_a = G_a[\gamma_a \lambda_a R_a I_F - \lambda'_a I]\check{f}_a, \qquad \bar{a} < a < 0,$$

*where $\bar{a}$ is defined in Lemma 2.9. These functions satisfy the uniform bound,*

$$\sup_{\bar{a} < a < 0} \|\check{f}'_a\|_V < \infty.$$

PROOF. The eigenfunction normalized with $\nu(\check{f}_a) = 1$ can be expressed as (28), where the potential kernel $G_a$ exists as a bounded linear operator from $L_\infty^W$ to $L_\infty^V$ by the previous two lemmas.

For any two values $a_1, a_2$ we then have

$$(31) \quad \check{f}_{a_2} - \check{f}_{a_1} = G_{a_2} s - G_{a_1} s = G_{a_2}[(\gamma_{a_1} - \gamma_{a_2})I - (R_{a_1} - R_{a_2})]G_{a_1} s.$$

Lemmas 2.9 and 2.8 imply that $\check{f}_{a_1} = G_{a_1} s$ is uniformly bounded. Convexity of $\Lambda$ implies that $|\gamma_{a_1} - \gamma_{a_2}|(a_2 - a_1)^{-1}$ is uniformly bounded for $\bar{a} \leq a_2 < a_1 < 0$, and it may then be verified using the mean value theorem that for some constant $b_0 < \infty$,

$$(a_2 - a_1)^{-1}(|\gamma_{a_1} - \gamma_{a_2}| + \|R_{a_1}\mathbf{1} - R_{a_2}\mathbf{1}\|_W) \leq b_0, \qquad \bar{a} \leq a_2 < a_1 < 0.$$

Applying Lemma 2.8 once more we conclude that

$$(a_2 - a_1)^{-1}\|\check{f}_{a_2} - \check{f}_{a_1}\|_V \leq b_1, \qquad \bar{a} \leq a_2 < a_1 < 0,$$

where $b_1 := b_0 \|\bar{G}W\|_V < \infty$. These bounds justify considering the limit $a_2 \downarrow a_1$ in (31) to obtain both the desired expression for $\check{f}'_a$ and the uniform bounds. □

LEMMA 2.11. *Suppose that the assumptions of Theorem 2.6 hold. Then $\check{f}_a \to 1$ and $\lambda_a \to 1$ as $a \uparrow 0$.*

PROOF. This follows from the uniform bound $G_a \leq \bar{R}$, and the formula for the limiting potential kernel,

$$G = \sum_{k=0}^{\infty} [R - s \otimes \nu]^k.$$

We have $\check{f}_a = G_a s \to Gs$ as $a \uparrow 0$, and it is known that $Gs \equiv 1$ since $\boldsymbol{\Phi}$ is Harris recurrent ([43], Theorem 5.1). □

PROOF OF THEOREM 2.6. Part (i) is given in Lemma 2.7; (ii) follows from Lemmas 2.8 and 2.9; and (iii) is given in Lemma 2.10.

To see that $\Lambda$ is smooth we argue as in [33, 34]: the g.p.e. $\gamma_a$ for $R_a$ is defined for $a \in (\bar{a}, 0)$ as the unique solution to

$$\nu[I\gamma - (R_a - s \otimes \nu)]^{-1} s = 1.$$

Hence smoothness follows from the inverse function theorem—see Proposition 4.8 of [33].

We now show that $\phi(a) = \Lambda'(a) \to \phi$ as $a \uparrow 0$. Based on Lemma 2.10 we have

$$h := \lim_{a \uparrow 0} \check{f}'_a = GX\check{f}_{0-},$$

with $X := (\gamma_{0-} \lambda_{0-})RI_F - \lambda'_{0-}I$. Moreover, $\nu(\check{f}'_a) = 0$ for each $a$, and hence from the uniform bounds on $\check{f}'_a$ combined with the dominated convergence theorem we have $\nu(h) = 0$.

Lemma 2.11 implies that $\gamma_{0-} = \lambda_{0-} = 1$ so that $X = [RI_F - \lambda'_{0-}I]$, and also $\check{f}_{0-} = 1$ so that

$$0 = \nu(h) = \nu GX\mathbf{1} = \mu(RF - \lambda'_{0-}),$$

where $\mu = \nu G$ is an unnormalized invariant measure. In particular $\mu R = \mu$, so that the expression above implies that $\mu(F) = \lambda'_{0-}\mu(\mathsf{X})$, or equivalently,

$$\phi := \pi(F) = \frac{\mu(F)}{\mu(\mathsf{X})} = \lambda'_{0-}. \qquad \square$$

The following result is a weak version of Varadhan's lemma [11].

PROPOSITION 2.12. *The following are equivalent for a nonnegative, measurable function $F : \mathsf{X} \to \mathbb{R}_+$, and any given initial condition $x_0 \in \mathsf{X}$:*

(i) *For some $c_0 > 0$,*

$$-\bar{I}(c_0) := \limsup_{n \to \infty} \frac{1}{n} \log \mathsf{P}_{x_0}\left\{\sum_{k=0}^{n-1} F(\Phi(k)) \geq nc_0\right\} < 0.$$

(ii) *For some $\theta_0 > 0$,*

$$\bar{\Lambda}(\theta_0) := \limsup_{n \to \infty} \frac{1}{n} \log \mathsf{E}_{x_0}\left[\exp\left(\sum_{k=0}^{n-1} \theta_0 F(\Phi(k))\right)\right] < \infty.$$

PROOF. The implication (ii) $\Rightarrow$ (i) is simply Chernoff's bound. Conversely, if (i) holds, then there exists $K_0 < \infty$ such that

$$\mathsf{P}_{x_0}\left\{\sum_{k=0}^{n-1} F(\Phi(k)) \geq nc_0\right\} \leq K_0 e^{-\bar{I}(c_0)n}, \qquad n \geq 0.$$

Consequently, since $F$ is assumed nonnegative-valued, we have for any $r > 1$,

$$\mathsf{P}_{x_0}\left\{\sum_{k=0}^{n-1} F(\Phi(k)) \geq nrc_0\right\} \leq \mathsf{P}_{x_0}\left\{\sum_{k=0}^{\lfloor nr \rfloor - 1} F(\Phi(k)) \geq nrc_0\right\}$$

$$\leq K_0 e^{-\bar{I}(c_0)(nr-1)}, \qquad n \geq 0.$$

Fix $\theta_0 < c_0^{-1}\bar{I}(c_0)$. On multiplying each side of this bound by $\exp(\theta_0 c_0 nr)$ we obtain

$$\mathsf{E}_{x_0}\left[\mathbb{I}\left\{\sum_{k=0}^{n-1} F(\Phi(k)) \geq nrc_0\right\} e^{\theta_0 c_0 n}\right] \leq K_1 e^{(\theta_0 c_0 - \bar{I}(c_0))nr}, \qquad n \geq 0,$$

where $K_1 := K_0 e^{\bar{I}(c_0)}$. Integrating both sides from $r = 1$ to $\infty$ we arrive at the bound, for each $n \geq 0$ and $\theta_0 < c_0^{-1} \bar{I}(c_0)$,

$$\mathsf{E}_{x_0}\left[\exp\left(\sum_{k=0}^{n-1} \theta_0 F(\Phi(k))\right)\right] \leq e^{c_0 \theta_0 n}\left(1 + K_1 \frac{c_0 \theta_0}{\bar{I}(c_0) - c_0 \theta_0} e^{-\bar{I}(c_0) n}\right).$$

This implies (ii) with $\bar{\Lambda}(\theta_0) \leq c_0 \theta_0$. $\square$

2.3. *Proof of Theorem* 2.2. Under the assumptions of Theorem 2.6 the function $\Lambda(a)$ is convex on $(\bar{a}, 0)$. Define the parameters,

$$\overline{c}_0 = \lim_{a \downarrow \bar{a}} \frac{d}{da} \Lambda(a); \qquad \overline{c}_1 = \lim_{a \uparrow 0} \frac{d}{da} \Lambda(a).$$

Theorem 2.6(iv) implies that $\overline{c}_1 = \phi$. By convexity we have $\overline{c}_0 \leq \overline{c}_1$, and this inequality is strict if $F$ is nondegenerate since then

$$(32) \qquad \sigma_a^2 := \frac{d^2}{da^2} \Lambda(a) > 0, \qquad a \in (\bar{a}, 0).$$

Equation (32) follows from [33], Proposition 2.4 (see also [34], Lemma 4.12).

For $c \in (\overline{c}_0, \overline{c}_1)$ the convex dual of $\Lambda$ is expressed

$$I(c) = \max_{a \in (\bar{a}, 0)}[ca - \Lambda(a)] = ca^* - \Lambda(a^*),$$

where $a^*$ is chosen so that $\frac{d}{da}\Lambda(a) = c$. The function $I$ serves as an LDP rate function within this range.

The LDP (3) for $c \in (\overline{c}_0, \phi)$ then follows from the multiplicative ergodic Theorem 2.3 and standard arguments [4, 11].

The proof of the exact LDP (14) is identical to that of the corresponding results in [33, 34], where

$$g_c = \frac{1}{a^* \sigma_{a^*}} \check{f}_{a^*},$$

$a^* \in (\bar{a}, 0)$ is again chosen so that $\frac{d}{da}\Lambda(a) = c$, $\check{f}_{a^*}$ is the eigenfunction satisfying the normalization (19) and $\sigma_{a^*}^2$ is defined in (32). The proof amounts to verification of the assumptions of [4], based on the multiplicative ergodic Theorem 2.3.

**3. Application to control-variates.** In this section we show how the method of control-variates can be used to construct a simulator that satisfies an LDP for the upper and lower tails even when the assumptions on $F$ in Theorem 2.2 are violated.

We begin with a general application of Theorem 2.2.

3.1. *Control-variates based on a Lyapunov function.* Suppose that (V3) holds, and consider the function $H := V - PV$. If $\pi(V) < \infty$, then invariance of $\pi$ implies that $\pi(H) = 0$, and hence the function $H$ can be used to construct a control-variate as described in the Introduction.

Recall that the assumption $\pi(V) < \infty$ means that the chain is $W$-regular of degree 2. In the following result we prefer to avoid this restriction and simply assume directly that $\pi(|H|) < \infty$ and that $\pi(H) = 0$. Under this assumption, with $H := V - PV$, define the sequence,

$$\Delta_n := n^{-1} \sum_{k=0}^{n-1} H(\Phi(k)), \qquad n \geq 0. \tag{33}$$

Positive Harris recurrence of $\Phi$ implies that $\Delta_n \to 0$ as $n \to \infty$ with probability 1 ([38], Theorem 17.0.1). The control-variate for simulation of a function $F$ based on the control-variate $H$ and a given parameter $\theta \in \mathbb{R}$ is given by

$$L_n(F_\theta) := L_n(F) - \theta L_n(H) = L_n(F) - \theta \Delta_n, \qquad n \geq 1.$$

In Theorem 3.1 we fix a function $F \in L_\infty^W$ together with constants $\theta_-, \theta_+$ each strictly greater than $\|F\|_W$. Define

$$F_- = F - \theta_- H, \qquad F_+ = F + \theta_+ H,$$

and the pair of estimators

$$\phi_n^- = L_n(F_-), \qquad \phi_n^+ = L_n(F_+), \qquad n \geq 1. \tag{34}$$

We denote by $\Lambda_-(a), \Lambda_+(a)$ the logarithm of the generalized principal eigenvalue for each of the scaled kernels

$$P_a^- = I_{e^{aF_-}} P, \qquad P_a^+ = I_{e^{aF_+}} P, \qquad a \in \mathbb{R}.$$

The convex duals of the functions $\{\Lambda_-, \Lambda_+\}$ are denoted $\{I_-, I_+\}$.

THEOREM 3.1. *Suppose that* (V3) *holds with* $W : \mathsf{X} \to [1, \infty)$ *near-monotone and* $V : \mathsf{X} \to (0, \infty)$ *everywhere finite, and suppose that* $\pi(H) = 0$, *where* $H := V - PV$. *Then, for a given function* $F \in L_\infty^W$, *there exists* $\varepsilon_0 > 0$ *such that:*

(i) *The lower LDP limit holds using* $\{\phi_n^+\}$:

$$\lim_{n \to \infty} n^{-1} \log \mathsf{P}\{\phi_n^+ \leq c\} = -I_+(c), \qquad c \in (\phi - \varepsilon_0, \phi).$$

(ii) *The upper LDP limit holds using* $\{\phi_n^-\}$:

$$\lim_{n \to \infty} n^{-1} \log \mathsf{P}\{\phi_n^- \geq c\} = -I_-(c), \qquad c \in (\phi, \phi + \varepsilon_0).$$

(iii) *If in addition $F$ is strongly nonlattice, then* (i) *and* (ii) *can be strengthened to the corresponding exact LDP limit analogous to* (14).

Parts (i) and (ii) of Theorem 3.1 combined imply that

$$\lim_{n\to\infty} n^{-1}\log \mathsf{P}_x\{[\phi_n - \phi] \in [-\varepsilon - \theta_+\Delta_n, \varepsilon + \theta_-\Delta_n]^c\}$$
$$= -\min(I_+(\phi - \varepsilon), I_-(\phi + \varepsilon)) < 0, \qquad 0 < \varepsilon < \varepsilon_0.$$

PROOF OF THEOREM 3.1. By normalization we can assume without loss of generality that $\pi(F) = 0$.

Define $W' = H + b\mathbb{I}_C$ so that by definition of $H$ we have $PV = V - W' + b\mathbb{I}_C$, and by (V3) we also have the lower bound,

$$W' = (V - PV) + b\mathbb{I}_C \geq (W - b\mathbb{I}_C) + b\mathbb{I}_C = W.$$

Both $F$ and $H$ belong to the function space $L_\infty^{W'}$.

We can write

$$F_- = F - \theta_- H = F - \theta_-(W' - b\mathbb{I}_C),$$

which implies that $-F_-$ is near-monotone whenever $\theta_- > \|F\|_W$ since we have the explicit bound

$$-F_- \geq (\theta_- - \|F\|_W)W - \theta_- b\mathbb{I}_C.$$

Similarly, $F_+$ is near-monotone whenever $\theta_+ > \|F\|_W$ since we can obtain the similar lower bound

$$F_+ = F + \theta_+(W' - b\mathbb{I}_C)$$
$$\geq -\|F\|_W W + \theta_+(W' - b\mathbb{I}_C) \geq (\theta_+ - \|F\|_W)W - \theta_+ b\mathbb{I}_C.$$

Moreover, in either case (13) holds for some $r > 0$, $r_0 < \infty$, since $W \in L_\infty^V$.

Hence the conclusions of Theorem 3.1 follow from Theorem 2.2. □

3.2. *Application to simulation of queues.* We now return to the reflected random walk (10) to illustrate the conclusions of Theorem 2.2.

The assumptions of Proposition 3.2 will be imposed throughout this section. We do not assume that $\mathsf{E}[|D(k)|^p] < \infty$ for any $p > 2$, so the CLT may not hold (see [23], [1], Theorem 3.2, Chapter 5 and [38], Chapter 17). Moreover, $\boldsymbol{\Phi}$ may not be geometrically ergodic since we do not assume that the distribution of $D(k)$ has exponential tails.

PROPOSITION 3.2. *Consider the reflected random walk* (10) *satisfying*

$$\delta := -\mathsf{E}[D(k)] > 0, \qquad \sigma_D^2 = \mathrm{Var}(D) < \infty, \qquad \mathsf{P}\{D(k) > 0\} > 0.$$

*Then:*

SIMULATING LARGE FUNCTIONS 23(i) *The following identity holds:*

$$PV(x) = V(x) - x + R(x), \qquad x \in \mathbb{R}_+, \tag{35}$$

where $V(x) = 1 + \frac{1}{2}\delta^{-1}(x^2 + \delta x)$ and $R$ is bounded. Hence (V3) holds with $W(x) = 1 + \frac{1}{2}x$, $x \in \mathbb{R}_+$.

(ii) *With $\tau_0$ equal to the first return time to the origin, we have for each $x \in \mathbb{R}_+$,*

$$\lim_{r \to \infty} r^{-1} \mathsf{E}_{rx}[\tau_0] = \delta^{-1} x, \qquad \lim_{r \to \infty} r^{-2} \mathsf{E}_{rx}\left[\sum_{k=0}^{\tau_0 - 1} \Phi(k)\right] = \tfrac{1}{2}\delta^{-1} x^2.$$

(iii) *A unique steady-state distribution $\pi$ exists satisfying $\int e^{\beta x} \pi(dx) = \infty$ for all $\beta > 0$ sufficiently large.*

(iv) *Let $\lambda_a$ denote the g.p.e. for the kernel $P_{f_a}$ with $f_a = e^{aF}$. Then $\lambda_a = \infty$ for all $a > 0$ when $F(x) \equiv x$.*

The proof of Proposition 3.2 is postponed to the end of this section. This result combined with Theorem 2.2 implies the LDP for the $M/M/1$ queue:

PROOF OF PROPOSITION 1.1. Either of the functions $F(x) = x$ or $F(x) = e^{\beta x}$ is near-monotone, with $\phi = \pi(F) < \infty$. Moreover, for $F(x) = x$ condition (V3) holds by Proposition 3.2(i), and with $F(x) = e^{\beta x}$ for a fixed $\beta \in (0, |\log(\rho)|)$, the functions $V = kF$, $W = F$ solve (V3) with $k = [1 - (\alpha e^\beta + (1-\alpha)e^{-\beta})]^{-1}$. Hence the one-sided LDP follows from Theorem 2.2. The proof of the LDP for positive $a$, with rate function satisfying $I(a) = 0$ for $a \geq \phi$, follows from Proposition 3.2(iv) combined with Proposition 2.12. □

Although Proposition 1.1 is stated for the $M/M/1$ queue, analogous conclusions hold for the general reflected random walk on $\mathbb{R}_+$ under the assumptions of Proposition 3.2. Part (i) asserts that (V3) holds, and hence the assumptions of Theorem 3.1 hold with $F(x) \equiv x$, so that the standard estimator satisfies a one-sided LDP.

We now show how the Lyapunov function $V$ can be used to construct a control-variate to obtain both upper and lower error bounds. Note that we do not know if $\pi(V) < \infty$ since we have not assumed that $D$ possesses a third moment. Consequently, we must use some other means to establish that $\pi(H) = 0$.

Under (V3) it follows from Proposition 2.1(i) that there exists a solution $\widehat{F} \in L_\infty^V$ to Poisson's equation (9). Moreover, it is known that the following scaling property holds:

$$\lim_{r \to \infty} r^{-2} \widehat{F}(rx) = J(x), \qquad x \in \mathbb{R}_+, \tag{36}$$

where $J$ is the fluid value function, $J(x) = \frac{1}{2}\delta^{-1}x^2$. The function $\widehat{F}$ is convex, unique up to an additive constant and can be chosen so that $\widehat{F}: \mathsf{X} \to \mathbb{R}_+$. The limit result (36) follows from Proposition 5.3 of [6] (see also [35], Theorem 16). Convexity is established in [5, 37] for network models.

Iterating Poisson's equation gives

$$P^n\widehat{F} = \widehat{F} + n\phi - \sum_{k=0}^{n-1} P^k F, \qquad n \geq 1.$$

It follows from the $f$-norm ergodic theorem [38] that for each initial condition $x \in \mathbb{R}_+$,

$$n^{-1}P^n\widehat{F}(x) = n^{-1}\widehat{F}(x) + \phi - n^{-1}\sum_{k=0}^{n-1} P^k F(x) \to 0, \qquad n \to \infty.$$

The quadratic growth and positivity of $\widehat{F}$ imply that we can find $\varepsilon > 0$ such that $\widehat{F}(x) \geq \varepsilon V(x) - 1$ for all $x$. Since $n^{-1}P^n\widehat{F}(x) \to 0$ as $n \to \infty$, we conclude that also $n^{-1}P^n V(x) \to 0$ as $n \to \infty$ for each $x$.

On setting $H = V - PV$ we see from (35) that the function $H$ can be written

(37) $$H(x) = x - R(x), \qquad x \in \mathbb{R}_+,$$

where $R: \mathsf{X} \to \mathbb{R}$ is bounded. In particular, it has linear growth so that $\pi(|H|) < \infty$. Moreover, we have

$$P^n V = V - \sum_{k=0}^{n-1} P^k H, \qquad n \geq 1,$$

and since $n^{-1}P^n V \to 0$ pointwise, we conclude that $\pi(H) = 0$. This justifies consideration of $F_\theta = F - \theta H$ in an asymptotically unbiased estimator of $\phi$.

It also follows from the representation (37) that for any given $r \geq 1$, the set $C_{F_\theta}(r)$ is compact whenever $\theta < 1$, and $C_{-F_\theta}(r)$ is compact whenever $\theta > 1$. This structure allows the application of Theorem 2.2 to obtain confidence bounds:

PROPOSITION 3.3. *Consider the reflected random walk* (10) *with* $\mathsf{E}[D(k)^2] < \infty$; $\mathsf{P}\{D(k) > 0\} > 0$; *and* $\delta > 0$. *Let* $F(x) = x$ *for* $x \in \mathbb{Z}_+$, *fix two parameters* $\theta_+ < 1$ *and* $\theta_- > 1$, *and define the pair of estimators* $\{\phi_n^-, \phi_n^+\}$ *via* (34) *with*

$$F_- = F - \theta_- H, \qquad F_+ = F - \theta_+ H.$$

*Then, there exists a pair of convex functions* $\{I_-, I_+\}$ *on* $\mathbb{R}$, *and a constant* $\varepsilon_0 > 0$ *such that for each initial condition* $x \in \mathsf{X}$,

$$\lim_{n\to\infty} n^{-1}\log \mathsf{P}_x\{\phi_n^+ \leq c\} = -I_+(c) < 0, \qquad c \in (\phi - \varepsilon_0, \phi),$$

$$\lim_{n\to\infty} n^{-1}\log \mathsf{P}_x\{\phi_n^- \geq c\} = -I_-(c) < 0, \qquad c \in (\phi, \phi + \varepsilon_0).$$

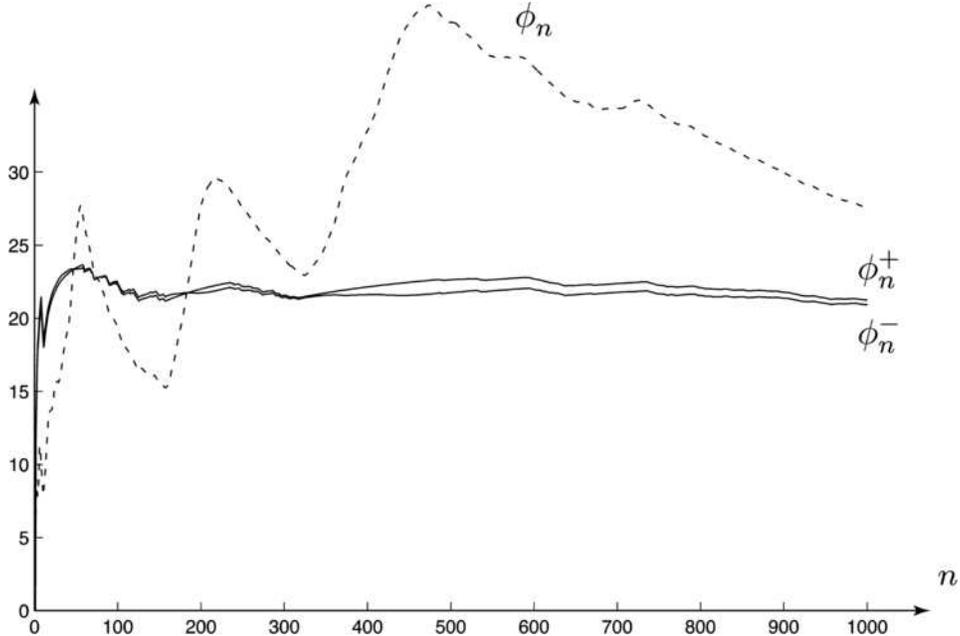

Fig. 2. *Monte Carlo estimates of $\phi := \pi(F)$ with $F(x) = x$ for $x \in \mathbb{R}_+$. The stochastic process $\boldsymbol{\Phi}$ is reflected random walk (10) with $\delta = -\mathsf{E}[D(k)] = 1$, and $\sigma_D^2 = 25$. The uncontrolled estimator exhibits large fluctuations around its steady-state mean. The upper and lower controlled estimators show less variability, and the bound $\phi_n^- < \phi_n^+$ is maintained throughout the run.*

Consequently, with $\{\Delta_n\}$ defined in (33), the following limit holds for each $\varepsilon \in (0, \varepsilon_0)$:
$$\lim_{n \to \infty} n^{-1} \log \mathsf{P}_x\{[\phi_n - \phi] \in [-\varepsilon + \theta_+ \Delta_n, \varepsilon + \theta_- \Delta_n]^c\}$$
$$= -\min(I_+(\phi - \varepsilon), I_-(\phi + \varepsilon)) < 0.$$

PROOF. The assumption that $\mathsf{P}\{D(k) > 0\} > 0$ is used to deduce that $\pi$ has support outside of the origin. For $\varepsilon > 0$ sufficiently small the set $C = [0, \varepsilon]$ is small, and satisfies the one-step minorization condition: for some $\delta > 0$, $P(x, \cdot) \geq \delta \nu(\cdot)$ with $\nu$ the point-mass at the origin. It follows that the functions $\{F, F_-, F_+\}$ are each nondegenerate.

The conclusions then follow from Theorem 3.1. □

An illustration of these controlled estimators is provided in Figure 2. The sequence $\mathbf{D}$ was chosen of the form $D(k) = A(k) - S(k)$, where $\mathbf{A}$ and $\mathbf{S}$ are mutually independent, i.i.d. sequences. Given nonnegative parameters $\mu, \alpha, \kappa$ we set
$$\mathsf{P}\{S(k) = (1 + \kappa)\mu\} = 1 - \mathsf{P}\{S(k) = 0\} = (1 + \kappa)^{-1},$$

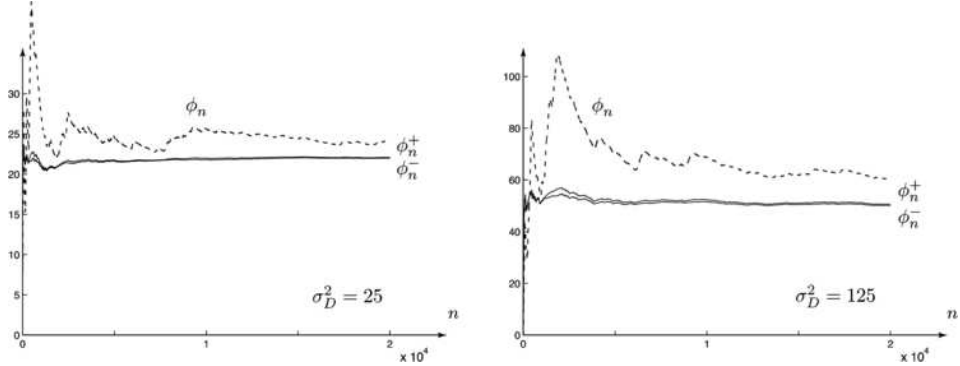

Fig. 3. *The plot at left shows the same simulation as shown in Figure 2, with the time horizon increased to $T = 20{,}000$. The plot at right shows the two controlled estimators along with the uncontrolled estimator when the variance is increased to $\sigma_D^2 = 125$. In each case the estimates obtained from the standard Monte Carlo estimator are significantly larger than those obtained using the controlled estimator, and the bound $\phi_n^- < \phi_n^+$ again holds for all large $n$.*

$$\mathsf{P}\{A(k) = (1+\kappa)\alpha\} = 1 - \mathsf{P}\{A(k) = 0\} = (1+\kappa)^{-1}.$$

Consequently, we have $\mathsf{E}[D(k)] = \mathsf{E}[A(k)] - \mathsf{E}[S(k)] = -(\mu - \alpha)$, and $\sigma_D^2 = \sigma_A^2 + \sigma_S^2 = (\mu^2 + \alpha^2)\kappa$. The simulation results shown in Figure 2 used $\mu = 4$, $\alpha = 3$ and $\kappa = 2$, so that $\delta = 1$ and $\sigma_D^2 = 25$.

The control-variate parameter values $\theta_- = 1.05$ and $\theta_+ = 1$ were used in the construction of $\{\phi_n^-, \phi_n^+\}$. While this value of $\theta_+$ violates the strict inequality $\theta_+ > 1$ required in Proposition 3.3, we have in this case

$$F_+(x) = x - \theta_+(x - R(x)) = R(x), x \in \mathbb{R}_+.$$

The function $R$ has mean zero and satisfies (13) when $\mathbf{D}$ has bounded support [or just a $(2+\epsilon)$-moment], so Theorem 2.2 implies that the lower LDP does hold using $\{\phi_n^+\}$ when $\theta_+ = 1$.

The plot at left in Figure 3 illustrates the simulation shown previously in Figure 2, with the time horizon increased to $T = 20{,}000$. The plot at right shows the controlled and uncontrolled estimators with $\kappa = 5$, and hence $\sigma_D^2 = 125$. The bounds $\phi_n^- < \phi_n^+ < \phi_n$ hold for all large $n$ even though all three estimators are asymptotically unbiased.

PROOF OF PROPOSITION 3.2. The function $R$ has the explicit form,

$$R(x) = \tfrac{1}{2}\delta^{-1}\sigma_D^2 - \tfrac{1}{2}\delta^{-1}\mathsf{E}[\{(x+D(k))^2 + \delta(x+D(k))\}\mathbb{I}(x < -D(k))],$$
$$x \in \mathsf{X}.$$

Under the second-moment assumption on $D(k)$ the function $R$ is bounded, since by Chebyshev's inequality,

$$\mathsf{E}[x^2 \mathbb{I}(x < -D(k))] \leq x^2 \mathsf{P}\{|D(k)| \geq x\} \leq \mathsf{E}[D(k)^2], \qquad x \in \mathbb{R}_+.$$

The second limit in (ii) follows from Proposition 5.3 of [6] [it can also be proved using (35) combined with the comparison theorem]. The proof of the first limit is similar and is omitted (see [7] for similar results).

The existence of $\pi$ satisfying $\pi(F) < \infty$ with $F(x) \equiv x$ follows from (35) and Proposition 2.1. On the interior of the set of $\beta \in \mathbb{R}$ satisfying $\pi(e^{\beta F}) < \infty$ we have by stationarity,

$$\pi(e^{\beta F}) = \mathsf{E}_\pi[\exp([\Phi(k) + D(k+1)]_+)] \geq \mathsf{E}_\pi[\exp(\Phi(k) + D(k+1))].$$

Hence by independence of $\Phi(k)$ and $D(k+1)$, the log moment generating functions for $\mathbf{\Phi}$ and $\mathbf{D}$ satisfy the bound

$$M(\beta) := \log(\pi(e^{\beta F})) \geq M(\beta) + M_D(\beta).$$

It follows that $M(\beta) = \infty$ when $M_D(\beta) > 0$, and this holds for large enough $\beta > 0$ under the assumption that $\mathsf{P}\{D(k) > 0\}$.

We now prove (iv). Suppose that in fact $\lambda_{\theta_0} < \infty$ for some positive $\theta_0$. It then follows that for $\Lambda_0 > \Lambda(\theta_0) := \log(\lambda_{\theta_0})$,

$$\sum_{n=0}^\infty \mathsf{E}_0\left[\exp\left(\sum_{k=0}^{n-1} \theta_0 \Phi(k) - \Lambda_0\right) \mathbb{I}\{\Phi(n) = 0\}\right] < \infty.$$

Define, with $\tau_0$ equal to the first return time to the origin,

$$h(x) := \mathsf{E}_x\left[\exp\left(\sum_{k=0}^{\tau_0-1} \theta_0 \Phi(k) - \Lambda_0\right)\right], \qquad x \in \mathsf{X}.$$

Then, from the previous bound,

$$h(0) = \sum_{n=0}^\infty \mathsf{E}_0\left[\exp\left(\sum_{k=0}^{n-1} \theta_0 \Phi(k) - \Lambda_0\right) \mathbb{I}\{\tau_0 = i\}\right] < \infty.$$

We next demonstrate that $h$ must be $\pi$-integrable.

For any $x \in \mathbb{R}_+$ we have

$$h(x) = \exp(\theta_0 F(x) - \Lambda_0)\mathsf{E}_x[h(\Phi_1)\mathbb{I}_{\Phi_1 \neq 0} + \mathbb{I}_{\Phi_1 = 0}],$$

from which it follows that the following identity holds for a bounded function $b_0 : \mathbb{R}_+ \to \mathbb{R}$:

$$\exp(\theta_0 F(x) - \Lambda_0)Ph = \exp(b_0(x))h(x), \qquad x \in \mathbb{R}_+.$$

This is a version of the drift condition (DV3) of [33, 34], which is far stronger than the drift condition (V3) of [38]. The comparison theorem of [38] implies that $\pi(h) < \infty$.

Next we obtain a lower bound on $h$ using Jensen's inequality:

$$\log h(x) \geq \mathsf{E}_x\left[\sum_{k=0}^{\tau_0-1}[\theta_0 \Phi(k) - \Lambda_0]\right].$$

Applying Proposition 3.2(ii), we conclude that the right-hand side is bounded from below by a quadratic function of $x$, giving a bound of the form, for some constant $b < \infty$,

$$\log h(x) \geq \tfrac{1}{2}\delta^{-1}x^2 - b(x+1), \qquad x \in \mathbb{R}_+.$$

This bound combined with Proposition 3.2(iii) implies that $\pi(h) = \infty$, which is a contradiction. □

**4. Conclusions.** We have seen that it is possible to establish strong LDP asymptotics for unbounded functions even when the assumptions of Donsker and Varadhan [46, 47] or the weaker geometric ergodicity assumption are violated. We are currently developing worst-case bounds when the statistics of the process are only partially known [32, 44], and we are also searching for ways of identifying explicit bounds on the rate function.

We are eager to develop these simulation techniques to better understand the value of the application of multiple control-variates for improved confidence bounds. Applications to network models are also considered in current research.

DEPARTMENT OF ELECTRICAL
AND COMPUTER ENGINEERING
AND THE COORDINATED SCIENCES LABORATORY
UNIVERSITY OF ILLINOIS AT URBANA–CHAMPAIGN
URBANA, ILLINOIS 61801
USA
E-MAIL: meyn@uiuc.edu